# Demand Adaptive Multi-Objective Electric Taxi Fleet Dispatching with Carbon Emission Analysis

Yiwen Song
School of Computer Science and Technology,
Beijing University of Posts and Telecommunications
Beijing, China

Ningning Sun
Xiamen No.1 High School of Fujian
Xiamen, China

Huimiao Chen
Sparkzone Institute
Beijing, China
chenhm15@tsinghua.org.cn

*Abstract*—As a foreseeable future mode of transport with lower emissions and higher efficiencies, electric vehicles (EVs) have received worldwide attention. For convenient centralized management, taxis are considered as the fleet with electrification priority. In this work, we focus on the study on electric taxis (ETs) dispatching, with consideration of picking up customers and recharging, based on real-world traffic data of a large number of taxis in Beijing. First, the assumed ET charging stations are located using the K-mean method. Second, based on the station locations and the order demands, which are in form of origin-destination (OD) pairs and extracted from the trajectory data, a dispatching strategy as well as the simulation framework is developed with consideration of reducing customer waiting time, mitigating ET charging congestion, and balancing order number distribution among ETs. The proposed method models the ET charging behaviors temporally discretely from the aspects of charging demands and availability of chargers, and further incorporates a centralized and intelligent fleet dispatching platform, which is capable of handling taxi service requests and arranging ETs' recharging in real time. The methodology in this paper is readily applicable to dispatching of different types of EV fleet with similar dataset available. Among the method, we use queueing theory to model the EV charging station waiting phenomena and include this factor into dispatching platform. Carbon emission is also surveyed and analyzed.

## I. INTRODUCTION

Global warming is threatening countries throughout the world, and countries have set a goal of limiting global temperature rise to well below 2 ºC. CO2 emissions are widely recognized as a large contributor to this phenomenon. As a kind of greenhouse gases, CO2 traps solar energy and warms the earth. Even though it's still likely that the increase of global temperature will still exceed 2 ºC, limiting CO2 emissions will decrease this possibility [1]. One critical aspect of limiting CO2 emissions will be the transportation sector. The number of vehicles in China have increased rapidly, and road vehicles are becoming one of the major contributors of CO2 emissions [2]. The transportation sector accounts for roughly 27 percent of the U.S. total greenhouse gas (GHG) emissions [3]. Another downside of vehicles is the large amounts of pollutants they generate. In Beijing, the annual mean concentration of fine particulates (PM 2.5) far exceeds the World Health Organization (WHO) guideline of annual mean limits, and vehicle tailpipe emissions contribute 31.1% of the total sources for the pollutant PM 2.5 [4]. Air pollutants emitted by vehicles greatly harm people's health, causing asthma, high blood pressure, lung cancer, premature deaths, etc. [5]. These all suggest that the current transportation systems are unsustainable.

Interests in electric vehicles (EVs) have increased due to advances in battery technologies, rising prices of petroleum, and growing concern over environmental issues [6], and adopting EVs is becoming a favorable option to combat the problems aforementioned. EVs are expected to play an important role in reducing air pollution, greenhouse gases, as well as health risks [5], promising far-reaching impacts. EVs powered by the European electricity mix could achieve a 10% to 24% decrease in global warming potential (GWP) [7]. Governmental support also helps the promotion of EVs. EVs are being proposed in China as an option to address the increasing energy demand from on-road transportation [8]. In the U.S., motor vehicle emission standards are explicitly set in clean air legislation, and policies to reduce the use of vehicles are designed. By imposing high taxes on fuel and giving subsidies to urban mass transit and intercity rail travel, European countries also effectively reduce vehicle use [9]. Plus, governments have made pledges that would boost the development of EVs. In 2012, China published a national plan that set a target of 5 million EVs in total by the year 2020 [10]. In May 2016, Indian government pledged to make all vehicles in India electricity-powered by 2030 [3].

Carbon emission analysis can help optimize the dispatching strategies of electric taxi (ET) fleet, and can help better assess the extent of global warming mitigation by EVs. Adopting different ET dispatching strategies can generate different results, helping us optimize goals including customer waiting time reduction or carbon emissions reduction. Through analysis of taxi dispatching simulation results, we can select the most effective strategy according to differed situations. Reference [11] provides insights into a rule-based approach of taxi dispatching. The authors simulate taxi dispatching using real data acquired from a taxi dispatch company, integrate information (e.g. the customer's location, time, the number of passengers, and the state of vacancy) to dispatch taxi, and collect various statistics including taxi utilization and average total waiting time of customers. The environmental impacts of EVs are analyzed from various aspects. In an effort to demonstrate how reductions of GHG emissions could help limit global warming, authors of [1] utilize a comprehensive probabilistic analysis focusing on quantifying GHG emission budgets for the first half of 21st century. Reference [7] develops a transparent life cycle inventory of conventional and electric vehicles and applies it to assess the impact EVs have on decreasing GWP. References [2] and [9] both evaluate the emissions of $CO_2$ and pollutants of vehicles in China. The authors of [8] point out that EVs in China involve multiple environmental issues as they use electricity generated mainly from coal, and compare the emissions of EVs with traditional vehicles (TVs). The authors of [2] base their study on a literature review and measured data, evaluating emissions using vehicle kilometers travelled (VKT) and emission factors. Reference [5] reviews the effects of EVs on GHG emissions, air quality, and human health. Authors of [4] comprehensively assess the ETs in Beijing in terms of energy efficiency, energy consumption, etc.

Despite the already done work, a more comprehensive dispatching method for EVs especially EV fleet like ET needs to be studied. In this work, we focus on the study on ETs' dispatching, with consideration of picking up customers and recharging, based on real-world traffic data of a large number of taxis in Beijing. The dispatching strategy as well as the simulation framework is developed with consideration of reducing customer waiting time, mitigating ET charging congestion, and balancing order number distribution among ETs. The proposed method models the ET charging behaviors from the aspects of charging demands and availability of chargers, and further incorporates a centralized and intelligent fleet dispatching platform, which is capable of handling taxi service requests and arranging ETs' recharging in real time. In the method, we use queueing theory to model the ET waiting phenomena in charging stations and consider this factor in dispatching strategy. Carbon emission is also surveyed and analyzed.

The remainder of the paper is organized as follows. Section II describes data and the station deployment. Section III introduce the real-time dispatching method and Section IV shows rating system for dispatching ETs. Section V surveys the $CO_2$ Emissions and Section VI gives the case. Section VI concludes.

## II. DATA DESCRIPTION AND THE DEPLOYMENT OF CHARGING STATIONS

To better reach our goals of optimizing the method of dispatching taxis and analyzing carbon emissions, we use the data of around 40 thousand taxis in Beijing from May 1st to May 27th in 2016, collected by smartphones. Note that because battery range is limited, we assume that ET drivers all use e-hailing or ride-sourcing smartphone applications to search for customers instead of wasting battery range on cruising on streets. The dataset records the status, including location, time and whether a taxi is carrying passengers about every 30 seconds, sorted by time. Table I is an example of the data record. We approximate a trip's origin or destination using the nearest data point that is within three minutes from the trip's start or end time if the trip's origin or destination lacks the information of GPS coordinates, and trips that have a traveling time less than two minutes are marked invalid.

TABLE I. TAXI TRAVEL RECORD SAMPLE

| Taxi ID | Time | Longitude | Latitude | Load |
|---|---|---|---|---|
| 2851 | 20160502111411 | 116.326775 | 39.896811 | 1 |

Recharging ETs is an important step in simulation and analysis, so we investigate the deployment of ET charging infrastructure as the prerequisite for fleet dispatching. The distribution of trip origins and destinations is significant for deploying charging stations for taxis, because ETs are preferred to be recharged at an intermediate stop between two trips. We apply K-means to cluster the locations of the origins of trips and locate charging station to each cluster at its centroid.

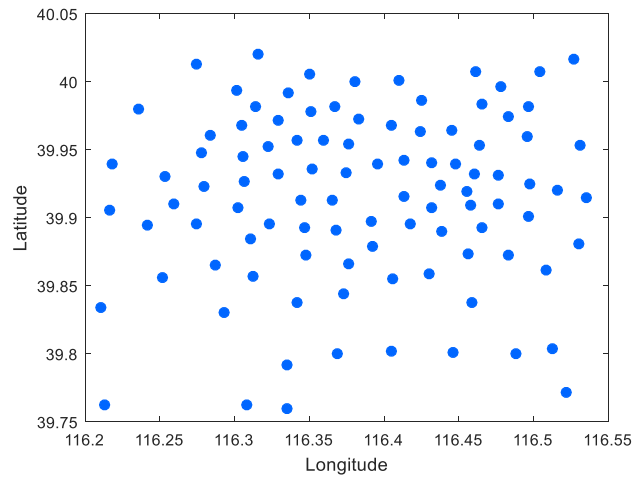

Fig. 1. Assumed Charging Station Locations.

## III. REAL-TIME ELECTRIC TAXI FLEET DISPATCHING

Based on the aforementioned dataset, we design a data-driven approach for dispatching ET fleet in this section. We conduct real-time simulation to gain insights on dispatching strategies in the context that both taxi drivers and passengers prefer to rely on the e-hailing or ride-sourcing smartphone applications to save battery range or for convenience. The centralized platform will match the ETs with customers based on the collected data on applications.

In the simulation, the taxi customers' travel demands are input to into the simulation dynamically. To facilitate the real-time operations of the taxi-dispatching platform and simplify the procedure, we divide the whole city into sub-regions according to the locations of charging stations (all geographic points sharing the same nearest charging stations form a sub-region). The sub-regions can also contribute to the identification of available chargers nearby for ETs with charging demand. Note that all the distances in this paper refer to the Manhattan distance because of the typically horizontal and vertical road layout in cities.

The simulation for real-time dispatching can be explained in the following steps:

*Step 1*: The centralized platform receives a new customer request and locates the sub-region that request belongs to. Then the platform selects all available taxis in the sub-region as candidates.

*Step 2*: Feed all candidates into reachability test, which identifies whether the candidates' battery states of charge (SOC) are enough to satisfy the customer's travel request.

*Step 3*: If no candidates pass the reachability test, go to Step 6. Otherwise, go to *Step 4*.

*Step 4*: Obtain all the candidates that passed the reachability test into the grading module, and calculate their scores. The grading module include pick-up distance, length of taxi's current empty time, driver's cumulative income level, and the matching degree between utilization level of the destination subregion and the taxi's SOC after the trip. Dispatch the candidate with the highest score to pick up the customer.

*Step 5*: After the taxi completes its task, check if its battery SOC is below a pre-predetermined threshold. If so, send the taxi to the charging station of the current sub-region and mark it as unavailable. Otherwise, mark the taxi as available, and the taxi will stay at the current location to wait for the next customer request. Note that if the charging station is full, the taxi needs to wait in the queue and is still marked as unavailable. Station operation chart record and update real-time operational states of charging stations.

*Step 6*: Add the demand to the waiting list which follows the first-in-first-out (FIFO) rule. Ride requests that have already been in the waiting list own higher priority than the newly coming requests.

*Step 7*: When the customer's waiting time exceeds a predetermined threshold, the platform will dispatch the taxis from adjacent areas. If there's an available taxi that can pass the reachability test, then dispatch the taxi to pick up the customer.

*Step 8*: When the customer's waiting time exceeds a predetermined cancelling threshold, decline the ride request.

IV. RATING SYSTEM FOR DISPATCHING ELECTRIC TAXIS

In this section, we provide the rating method for dispatching taxis to satisfy the travel demands. We take into account the vehicle-to-origin distance (approximately directly proportional to the passenger waiting time), the cumulative income of taxi drivers, and the match of destination charging station occupy and the ET SOC after a trip.

The designed formula for score calculation is given in (1), shown as below.

$$Score = f(demand)w_1 d + w_2 I + w_3 O \quad (1)$$

where, $d$ represents the distance from ET to the trip origin, $I$ represents the total income of the ET in a given time period, and $O$ denotes the match degree between destination charging station occupy and the SOC of the ET after a trip. $w_1$, $w_2$ and $w_3$ are weights. $f(demand)$ is an adjustable function to make the weight adaptive to the travel demand. The formula of $f$ is given in (2) In this paper, we summary a demand curve as the input of the function. By using this function, when the demand is high, the distance will be put at a more important place so that the fulfill rate can be guaranteed. Note that a score is for an ET and a trip.

In this way, for a given trip, the system can then calculate scores of all ETs and select the best one to finish the order.

$$f(demand) = x\sqrt{(demand)} \quad (2)$$

To quantify the matching degree $O$ we first calculate the mean waiting time of charging station, by leveraging the approximation for an M/G/s queueing system developed in [12] and [13], expressed as below:

$$W^{M/G/s}{}_1 = \frac{(1+\xi^2)W^{M/M/s}W^{M/D/s}}{2\xi^2 W^{M/D/s}+(1-\xi^2)W^{M/M/s}} \quad (3)$$

$$W^{M/G/s}{}_2 = \lim_{\xi \to \infty} W^{M/G/s}{}_1 = \frac{W^{M/M/s}W^{M/D/s}}{2W^{M/D/s}-W^{M/M/s}} \quad (4)$$

$$W^{M/M/s} = \frac{\lambda^s}{(s-1)!(s\mu-\lambda)^2\mu^{s-1}}\left[\sum_{Z=0}^{s-1}\frac{1}{z!}\left(\frac{\lambda}{\mu}\right)^2 + \frac{\lambda^s}{(s-1)!(s\mu-\lambda)\mu^{s-1}}\right]^{-1} \quad (5)$$

$$W^{M/D/s} = \frac{1}{2}\left[1 + H\frac{s\mu-\lambda}{\lambda}\left(1-e^{-\frac{\lambda(s-1)}{H(s\mu-\lambda)(s+1)}}\right)\right]W^{M/M/s} \quad (6)$$

$$H = \frac{s-1}{16s}\left[\left(\frac{10s+8}{2}\right)^{\frac{1}{2}} - 2\right] \quad (7)$$

$$\xi = \sigma\mu \quad (8)$$

In (3)-(8), $W^{M/G/s}$ denotes the mean waiting time of an M/G/s queueing system, and $W^{M/M/s}$ and $W^{M/D/s}$ are the mean waiting time of the corresponding M/M/s and M/D/s queueing systems, respectively; $\mu$ and $\sigma$ are the reciprocals of the mean charging time of ETs and the standard deviation of charging time, respectively. $\lambda$ denotes the arrival rate of ETs, and in order to make the stimulation more practical, the ET arrivals should follow Poisson distribution with a time-varying $\lambda$, i.e., $\lambda(t)$. As a result, we obtain $\lambda$ by collecting the data from charging situation in each station at intervals of 60 minutes and discretely regard the arrival rate within an hour as a constant. When the standard deviation of charging time equals 0, the formula of counting its reciprocal $\sigma$ is meaningless. Hence, we make the standard deviation tend to 0 and $\sigma$ tend to $\infty$, then compute the limit of the equation (2), as (3) shows. Consequently, when $\sigma$ equals to 0, we calculate the mean waiting time through equation (3). Otherwise, we adopt equation (2).

Then we transform the mean waiting time to the station utilization level and quantify the matching degree by (9) and (10).

$$u_n = \frac{W^{M/G/s}}{C_1+W^{M/G/s}} \quad (9)$$

$$O = C_2(u_n - o_e)^2 + C_3 \quad (10)$$

Since SOCs $\in [0,1]$, in order to compare SOCs and waiting time and match them, we map the average waiting time from [0, $+\infty$] to [0, 1] through formula (9) and denote it as station utilization level $u_n$. In (9), the coefficient $C_1$ should be valued according to the specific scenes and practical intends. In this paper, we denote $C_1$ as the mean value of $W^{M/G/s}$. Taxis with low SOC need to be charged as soon as possible, so they prefer stations with shorter waiting time. On the other hand, taxis with higher SOC are in no hurry to recharge their batteries, so there is no need of short waiting time. Consequently, smaller difference between SOC and waiting time level implies higher match between them and higher efficiency of charging stations. In (10), we take coefficients $C_2$ and $C_3$ based on the data to ensure that possibly higher matching degree is achieved. In such a way, it improves the possibility for low SOC ETs to get orders, guarantees that they can get charged as soon as possible when they finish the trip, and avoids ETs with high SOC taking up low SOC ETs' charging resources.

## V. SURVEY AND ANALYSIS OF $CO_2$ EMISSIONS

It's a common belief that EVs generate less carbon emissions than TVs which are powered by internal combustion engines, as TVs are powered by fossil fuels and produce a tremendous amount of green-house gases, e.g., carbon dioxide. One thing that needs to be considered is that the adoption of EVs does not mean zero emissions; the emissions of EVs depend on a nation's electricity generation structure. Though renewable energy continues to gain an increasing share in the overall energy source, fossil fuels including coal, oil and gas still account for the majority of energy consumption and the generation of electricity. To compare the carbon emissions of TVs and EVs, we summarize $CO_2$ emissions based on data of the constitution of China's electricity generation structure, as well as data of the generation of $CO_2$ emissions. The results show that TVs generate about twice as much carbon emissions as electric cars per mile, proving that electric cars are a more environmental-friendly transport mode than their fossil fuel-powered counterparts.

Taxis, when they are all fossil fuel powered, generate 10 times as much carbon emissions as private cars per day on average due to their traveled miles. Therefore, based on the results aforementioned, replacing traditional taxis (TTs) with ETs will contribute to a dramatic reduction of carbon emissions.

### A. China's Electricity Mix

The data of China's Electricity Mix in 2016 [14] are shown below in Table II. Based on the data of China's electricity mix, we can be able to determine how much carbon emissions EVs generate.

### B. $CO_2$ Emissions of TVs

We use three typical types of TVs, i.e., Toyota Yaris L, Chevrolet Sonic, Ford Escape [15], and use the average to represent the oil consumption and $CO_2$ emissions of the TV population. The data are shown in Table III.

### C. $CO_2$ Emissions of EVs

We use two typical types of EVs, i.e., 2016 Geely Emgrand ("Dihao") pure-electric sedan [17] and BYD Chin Pro DM [18], and use the average to represent the electricity consumption and $CO_2$ emissions of all electric cars. The data are shown in Table IV.

### D. Summary

From the data above, we can clearly see that TVs generate more than twice as much emissions as EVs do. Given the fact that taxis drive far more than private cars every day, turning TVs to ETs will reduce carbon emissions dramatically and therefore create a healthier urban environment.

TABLE II. CHINA'S ELECTRICITY MIX*

| Type | Thermal Power | Hydropower | Wind Power | Solar Power | Pumped Storage | Nuclear Power | Traditional Energy |
|---|---|---|---|---|---|---|---|
| Proportion | 72.1% | 18.86% | 3.95% | 1.08% | 0.5% | 3.49% | 72.1% |

*Carbon dioxide emissions of renewable energy are neglected here.

TABLE III. TANK CAPACITY, VEHICLE RANGE, OIL CONSUMPTION AND $CO_2$ EMISSIONS OF TVs

| Characteristics | Toyota Yaris L | Chevrolet Sonic | Ford Escape |
|---|---|---|---|
| Tank Capacity (L) | 43.91 | 46.18 | 59.43 |
| Vehicle Range (km) | 653.39 | 608.33 | 656.61 |
| Oil Consumption (L/100km) | 6.72 | 7.59 | 9.95 |
| Average Oil Consumption (L/100km) | 8.09 | | |
| CO2 Emissions* (kg/100km) | 18.61 | | |

*Calculated based on: 2.3 kg $CO_2$ / 1 L oil [16].

TABLE IV. TANK CAPACITY, VEHICLE RANGE, OIL CONSUMPTION AND $CO_2$ EMISSIONS OF EVs

| Characteristics | Toyota Yaris L | Chevrolet Sonic |
|---|---|---|
| Battery Capacity (kWh) | 45.3 | 14.38 |
| Distance Covered Per Trip* (km) | 286.7 | 82 |
| Electricity Consumption (kWh/100km) | 15.8 | 17.5 |
| Average Electricity Consumption (kWh/100 km) | 16.65 | |
| Average Coal Consumption** (kg/100km) | 3.73 | |

| | |
|---|---|
| $CO_2$ Emissions*** (kg/100km) | 8.88 |

*Each trip is fully electricity-powered; before each trip the battery is full and each trip uses up all the electricity storage.
**Calculated based on China's electricity mix and the assumption that all thermal power generation comes from coal.
***Calculated based on: 2.38 kg CO2 / 1kg coal [19], and 3.21 kWh / 1 kg coal [20].

## VI. CASE STUDIES AND SIMPLE ANALYSIS

Here, we conduct the sensitivity analyses with respect to $f(demand)$. Specifically, we run multiple simulations under various taxi fleet configurations, as shown in Table V. In this case, we allocate 9000 ETs in the area of about within the 5th ring road of the city of Beijing. The driving range of an ET is 250 km and the battery size is 38 kWh.

TABLE V. RESULTS UNDER DIFFERENT F(DEMAND)

| $x$ in $f(demand)$ | Passenger Mean Waiting Time (min) | Gini Coefficient | Fulfill Rate (%) | Mean Distance (km) |
|---|---|---|---|---|
| 0.001 | 10.6661 | 0.1923 | 0.9659 | 344.5267 |
| 0.01 | 11.2435 | 0.2038 | 0.9533 | 339.8923 |
| 0.1 | 10.2461 | 0.1734 | 0.9625 | 343.4098 |

From Table V, we can see that by using our dispatching strategy, the fulfill rate can be maintained at a high level, the gini coefficient of drivers are quite low and the waiting time of passengers are lower than 12 minutes. But we find just increasing the $x$ value in $f$ which means a higher weight of distance from the vehicle to the origin of a trip can not obviously improve the long-term fulfill rate and passengers' waiting time because in the simulation the number of ETs with close distance will diminish.

By using ETs, the $CO_2$ emission will be reduced about from 56,946,600 kg to 27,172,800 kg based on the data given in last section, which means a 52.3% decreasing.

## VII. CONSLUSIONS

In this work, we propose an ET rating system and further build an ET dispatching platform. The comprehensive dispatching strategy shows its advantage in the case simulation. Insights such as that excessive preference of a factor such as vehicle-to-origin distance may not improve the intuitively related indices are derived. In future work, we would like to study more accurate weights for the score calculation.